\documentclass{amsart}
\usepackage{amsmath,amssymb}
\usepackage{yhmath}

\newcommand{\bdis}{\begin{displaymath}}
\newcommand{\edis}{\end{displaymath}}
\newcommand{\be}{\begin{equation}}
\newcommand{\ee}{\end{equation}}
\newcommand{\mbb}{\mathbb}
\newcommand{\mcal}{\mathcal}

\newcommand{\vp}{\varphi}

\newcommand{\zf}{\zeta\left(\frac{1}{2}+it\right)} 
\newcommand{\zfu}{\zeta\left(\frac{1}{2}+iu\right)}

\newcommand{\bW}{\bar{W}}

\newcommand{\FR}{\frac{x^n+y^n}{z^n}} 
\newcommand{\FRI}{\frac{z^n}{x^n+y^n}}

\theoremstyle{definition}

\theoremstyle{remark}
\newtheorem{remark}[]{Remark}

\newtheorem*{mydef11}{{\bf Theorem 1}}

\newtheorem*{mydef12}{{\bf Theorem 2}}

\newtheorem*{mydef13}{{\bf Theorem 3}}

\newtheorem*{mydef14}{{\bf Theorem 4}}

\newtheorem*{mydef15}{{\bf Theorem 5}} 

\newtheorem*{mydef16}{{\bf Theorem 6}}

\newtheorem*{mydef51}{{\bf Lemma 1}}

\newtheorem*{mydef52}{{\bf Lemma 2}}

\newtheorem*{mydef53}{{\bf Lemma 3}}

\newtheorem*{mydef81}{{\bf Property 1}}

\newtheorem*{mydef82}{{\bf Property 2}}

\numberwithin{equation}{section}

\begin{document}

\title[Jacob's ladders, $\zeta$-transformation  \dots]{Jacob's ladders, $\zeta$-transformation of the Fourier orthogonal system (2014) and new infinite sets of $\zeta$-equivalents of the Fermat-Wiles theorem}

\author{Jan Moser}

\address{Department of Mathematical Analysis and Numerical Mathematics, Comenius University, Mlynska Dolina M105, 842 48 Bratislava, SLOVAKIA}

\email{jan.mozer@fmph.uniba.sk}


\begin{abstract}
In this paper we obtain new infinite sets of $\zeta$-equivalents of the Fermat-Wiles theorem based on the elementary Fourier orthogonal system, Riemann's zeta-function and Jacob's ladders. 
\end{abstract}
\maketitle

\section{Introduction} 

In this paper, that continues our series \cite{6} -- \cite{17}, we use our method of constructing new functionals on the orthogonal systems, especially on the classical Fourier system. We obtain new set of $\zeta$-equivalents of the Fermat-Wiles theorem. Simultaneously, we obtain a point of contact between Riemann's zeta-function, Hardy-Littlewood integral (1918), Jacob's ladders, elementary Fourier orthogonal system and Fermat-Wiles theorem.  

\subsection{} 

Let us remind our Theorem (see \cite{4}, (2.1) -- (2.4)): For every fixed $L_2$-orthogonal system 
\be \label{1.1} 
\{f_m(t)\}_{m=1}^\infty,\ t\in [0,2l],\ l=o\left(\frac{T}{\ln T}\right),\ T\to\infty , 
\ee 
there is the continuum sets of $L_2$-orthogonal systems 
\be \label{1.2} 
\begin{split}
& \{F_m(t;T,k,l)\}_{m=1}^\infty= \\ 
& \left\{f_m((\vp_1^k(t)-T))\prod_{r=0}^{k-1}|\tilde{Z}[\vp_1^r(t)]|\right\},\ t\in [\overset{k}{T},\overset{k}{\wideparen{T+2l}}]
\end{split}
\ee 
for every fixed $k\in\mbb{N}$, i.e. the following formula is valid 
\be \label{1.3} 
\begin{split}
& \int_{\overset{k}{T}}^{\overset{k}{\wideparen{T+2l}}}f_m((\vp_1^k(t)-T))f_n((\vp_1^k(t)-T))\prod_{r=0}^{k-1}|\tilde{Z}[\vp_1^r(t)]|^2{\rm d}t= \\ 
& \begin{cases}
	0,\ m\not=n, \\ A_m,\ m=n
\end{cases},\ A_m=\int_0^{2l}f_m^2(t){\rm d}t,\ A_m=A_m(l;f_m), 
\end{split} 
\ee

where 
\be \label{1.4} 
\tilde{Z}^2(t)=\frac{|\zf|^2}{\omega(t)},\ \omega(t)=\left\{1+\mcal{O}\left(\frac{\ln \ln t}{\ln t}\right)\right\}\ln t. 
\ee 

\begin{remark}
By this theorem it is true, that a new $L_2$-orthogonal system corresponds to every sufficiently big $T>0$. This new system is defined on the segment $[\overset{k}{T},\overset{k}{\wideparen{T+2l}}]$. Since 
\be \label{1.5} 
\overset{k}{T}=\vp_1^{-k}(T),\ \overset{k}{\wideparen{T+2l}}=\vp_1^{-k}(T+2l), 
\ee  
then we have the continuum set of new $L_2$-orthogonal systems generated by the Riemann's function $\zf$ and by the Jacob's ladders. 
\end{remark} 

\begin{remark}
For the purposes of this paper we shall assume that $l$ is independent on $T$, comp. (\ref{1.1}). 
\end{remark} 

\subsection{} 

We show in this paper that the formula\footnote{See (\ref{1.3})} connected with the normalization 
\be \label{1.6} 
\begin{split}
& \int_{\overset{k}{T}}^{\overset{k}{\wideparen{T+2l}}}f_m^2((\vp_1^k(t)-T))\prod_{r=0}^{k-1}|\frac{|\zeta(\frac 12+i\vp_1^r(t))|^2}{\omega[\vp_1^r(t)]}{\rm d}t=A_m, \\ 
& A_m=\int_0^{2l}f_m^2(t){\rm d}t,\ l>0,\ f_m(t)\in\{f_m(t)\}_{m=1}^\infty
\end{split}
\ee 
is good as the resource of new sets of $\zeta$-equivalents of the Fermat-Wiles theorem. 
 
\begin{remark}
For this purpose it is sufficient to assume that $[0,2l]\ni t\mapsto f_m(t)$ is continuous for every $m$, that is 
\be \label{1.7} 
f_m(t)\in \{f_m(t)\}_{m=1}^\infty\in \left\{\{f_m(t)\}_{m=1}^\infty\right\}_{C\perp [0,2l]}, 
\ee  
where the last symbol expresses the set of all continuous orthogonal systems on the segment $[0,2l]$. 
\end{remark} 

We obtain here, for example, the following result. Let the symbol 
\be \label{1.8} 
\FR,\ x,y,z,n\in\mbb{N},\ n\geq 3
\ee 
denote the Fermat's rationals. Then the following $\zeta$-condition 
\be \label{1.9} 
\begin{split}
& \lim_{\tau\to\infty}\frac{1}{\ln^k\tau}\int_{[W(\FR,\tau)]^k}^{[W(\FR,\tau)+2l]^k}f_m^2\left(\vp_1^k(t)-W\left(\FR,\tau\right)\right)\times \\ 
& \prod_{r=0}^{k-1}\left|\zeta\left(\frac 12+i\vp_1^r(t)\right)\right|^2{\rm d}t\not=1
\end{split}
\ee 
on the set of all Fermat's rationals expresses the next $\zeta$-equivalent of the Fermat-Wiles theorem for every fixed $k,m\in\mbb{N}$ and $l>0$\footnote{Comp. Remark 3.}, where we have used the following symbol 
\be \label{1.10} 
W(x,\tau)=W(x,\tau;k,m)=\exp\left\{\frac{x}{\sqrt[k]{A_m}}\ln\tau\right\} 
\ee  
with the substitution 
\bdis 
x\to \FR. 
\edis 

\subsection{} 

First, let us remind that to the classical Fourier orthogonal system 
\be \label{1.11} 
\left\{1,\cos\frac{\pi}{l}t, \sin\frac{\pi}{l}t, \dots , \cos\frac{\pi}{l}mt, \sin\frac{\pi}{l}mt, \dots\right\},\ t\in[0,2l] 
\ee 
corresponds the following set of orthogonal $\zeta$-transformations\footnote{See \cite{4}, (3.2).} 
\be \label{1.12} 
\begin{split}
& \left\{\prod_{r=0}^{k-1}|\frac{|\zeta(\frac 12+i\vp_1^r(t))|}{\sqrt{\omega[\vp_1^r(t)]}}, \dots, \right. \\ 
& \left. \left(\prod_{r=0}^{k-1}|\frac{|\zeta(\frac 12+i\vp_1^r(t))|}{\sqrt{\omega[\vp_1^r(t)]}}\right)\cos\left(\frac{\pi}{l}m(\vp_1^r(t)-T)\right), \right. \\ 
& \left. \left(\prod_{r=0}^{k-1}|\frac{|\zeta(\frac 12+i\vp_1^r(t))|}{\sqrt{\omega[\vp_1^r(t)]}}\right)\sin\left(\frac{\pi}{l}m(\vp_1^r(t)-T)\right), \dots
\right\}, \\ 
& t\in [\overset{k}{T},\overset{k}{\wideparen{T+2l}}],\ T\to\infty 
\end{split}
\ee 
for every fixed $k\in \mbb{N}$. 

Next, the followig set of normalization integrals corresponds\footnote{See (\ref{1.6}).} to the set (\ref{1.12}) 
\be \label{1.13} 
\int_{\overset{k}{T}}^{\overset{k}{\wideparen{T+2l}}}\prod_{r=0}^{k-1}\frac{|\zeta(\frac 12+i\vp_1^r(t))|^2}{\omega[\vp_1^r(t)]}{\rm d}t=2l, 
\ee  
\be \label{1.14} 
\int_{\overset{k}{T}}^{\overset{k}{\wideparen{T+2l}}}\prod_{r=0}^{k-1}\frac{|\zeta(\frac 12+i\vp_1^r(t))|^2}{\omega[\vp_1^r(t)]}\cos^2\left(\frac{\pi}{l}m(\vp_1^k(t)-T)\right){\rm d}t=l, 
\ee 
\be \label{1.15} 
\int_{\overset{k}{T}}^{\overset{k}{\wideparen{T+2l}}}\prod_{r=0}^{k-1}\frac{|\zeta(\frac 12+i\vp_1^r(t))|^2}{\omega[\vp_1^r(t)]}\sin^2\left(\frac{\pi}{l}m(\vp_1^k(t)-T)\right){\rm d}t=l, 
\ee 
where the last formula is a consequence of the first two. 

In this paper we obtain, for example, the following new $\zeta$-equivalents of the Fermat-Wiles theorem. 

\begin{itemize}
	\item[(A)] The following $\zeta$-equivalent 
	\be \label{1.16} 
	\lim_{T\to\infty}\frac{1}{\ln^kT}\int_{\overset{k}{T}}^{[T+2\FR]^k}\prod_{r=0}^{k-1}\left|\zeta\left(\frac 12+i\vp_1^r(t)\right)\right|^2{\rm d}t\not=2,
	\ee  
	corresponds to the integral (\ref{1.13}). 
	\item[(B)] The following $\zeta$-equivalent 
	\be \label{1.17} 
	\begin{split}
	& \lim_{T\to\infty}\frac{1}{\ln^kT}\int_{\overset{k}{T}}^{[T+2\FR]^k}\prod_{r=0}^{k-1}\left|\zeta\left(\frac 12+i\vp_1^r(t)\right)\right|^2\times \mcal{F}^2{\rm d}t\not=1, \\ 
	& \mcal{F}\in \left\{\cos\left(\pi m\FR(\vp_1^r(t)-T)\right), \sin\left(\pi m\FR(\vp_1^r(t)-T)\right)\right\}
	\end{split} 
	\ee  
	corresponds to the integrals (\ref{1.14}), (\ref{1.15}) for every fixed $k,m\in\mbb{N}$. 
\end{itemize} 

\subsection{} 

To the function, for example, 
\be \label{1.18} 
f_m(t)=\cos\left(\frac{\pi}{l}mt\right),\ m\in\mbb{N} 
\ee  
corresponds the following $\zeta$-equivalent\footnote{See (\ref{1.9}).} as a consequence of (\ref{1.6}), (\ref{1.18}), $A_m=l$ 
\be \label{1.19} 
\begin{split}
	& \lim_{\tau\to\infty}\frac{1}{\ln^k\tau}\int_{[\bW(\FR,\tau)]^k}^{[\bW(\FR,\tau)+2l]^k}\left(\prod_{r=0}^{k-1}\left|\zeta\left(\frac 12+i\vp_1^r(t)\right)\right|^2\right)\times \\ 
	& \cos^2\left(\frac{\pi}{l}m\left(\vp_1^k(t)-\bW\left(\FR,\tau\right)\right)\right) {\rm d}t\not=1
\end{split}
\ee 
where\footnote{Comp. (\ref{1.10}).} 
\be \label{1.20} 
\bW\left(x,\tau\right)=\exp\left[\frac{x}{\sqrt[k]{l}}\ln\tau\right] . 
\ee  

\begin{remark}
Consequently, we have two different $\zeta$-equivalents, namely the formula (\ref{1.17}) and the formula (\ref{1.19}), that we have constructed for the function (\ref{1.18}). These two $\zeta$-equivalents constitute the first dual of $\zeta$-equivalents, that is the infinite set of such duals, since $m\in\mbb{N}$. 
\end{remark}

\section{Jacob's ladders: notions and basic geometrical properties}  

\subsection{}

In this paper we use the following notions of our works \cite{5} -- \cite{9}: 
\begin{itemize}
\item[{\tt (a)}] Jacob's ladder $\vp_1(T)$, 
\item[{\tt (b)}] direct iterations of Jacob's ladders 
\bdis 
\begin{split}
	& \vp_1^0(t)=t,\ \vp_1^1(t)=\vp_1(t),\ \vp_1^2(t)=\vp_1(\vp_1(t)),\dots , \\ 
	& \vp_1^k(t)=\vp_1(\vp_1^{k-1}(t))
\end{split}
\edis 
for every fixed natural number $k$, 
\item[{\tt (c)}] reverse iterations of Jacob's ladders 
\be \label{2.1}  
\begin{split}
	& \vp_1^{-1}(T)=\overset{1}{T},\ \vp_1^{-2}(T)=\vp_1^{-1}(\overset{1}{T})=\overset{2}{T},\dots, \\ 
	& \vp_1^{-r}(T)=\vp_1^{-1}(\overset{r-1}{T})=\overset{r}{T},\ r=1,\dots,k, 
\end{split} 
\ee   
where, for example, 
\be \label{2.2} 
\vp_1(\overset{r}{T})=\overset{r-1}{T}
\ee  
for every fixed $k\in\mbb{N}$ and every sufficiently big $T>0$. We also use the properties of the reverse iterations listed below.  
\be \label{2.3}
\overset{r}{T}-\overset{r-1}{T}\sim(1-c)\pi(\overset{r}{T});\ \pi(\overset{r}{T})\sim\frac{\overset{r}{T}}{\ln \overset{r}{T}},\ r=1,\dots,k,\ T\to\infty,  
\ee 
\be \label{2.4} 
\overset{0}{T}=T<\overset{1}{T}(T)<\overset{2}{T}(T)<\dots<\overset{k}{T}(T), 
\ee 
and 
\be \label{2.5} 
T\sim \overset{1}{T}\sim \overset{2}{T}\sim \dots\sim \overset{k}{T},\ T\to\infty.   
\ee  
\end{itemize} 

\begin{remark}
	The asymptotic behaviour of the points 
	\bdis 
	\{T,\overset{1}{T},\dots,\overset{k}{T}\}
	\edis  
	is as follows: at $T\to\infty$ these points recede unboundedly each from other and all together are receding to infinity. Hence, the set of these points behaves at $T\to\infty$ as one-dimensional Friedmann-Hubble expanding Universe. 
\end{remark}  

\subsection{} 

Let us remind that we have proved\footnote{See \cite{9}, (3.4).} the existence of almost linear increments 
\be \label{2.6} 
\begin{split}
& \int_{\overset{r-1}{T}}^{\overset{r}{T}}\left|\zf\right|^2{\rm d}t\sim (1-c)\overset{r-1}{T}, \\ 
& r=1,\dots,k,\ T\to\infty,\ \overset{r}{T}=\overset{r}{T}(T)=\vp_1^{-r}(T)
\end{split} 
\ee 
for the Hardy-Littlewood integral (1918) 
\be \label{2.7} 
J(T)=\int_0^T\left|\zf\right|^2{\rm d}t. 
\ee  

For completeness, we give here some basic geometrical properties related to Jacob's ladders. These are generated by the sequence 
\be \label{2.8} 
T\to \left\{\overset{r}{T}(T)\right\}_{r=1}^k
\ee 
of reverse iterations of the Jacob's ladders for every sufficiently big $T>0$ and every fixed $k\in\mbb{N}$. 

\begin{mydef81}
The sequence (\ref{2.8}) defines a partition of the segment $[T,\overset{k}{T}]$ as follows 
\be \label{2.9} 
|[T,\overset{k}{T}]|=\sum_{r=1}^k|[\overset{r-1}{T},\overset{r}{T}]|
\ee 
on the asymptotically equidistant parts 
\be \label{2.10} 
\begin{split}
& \overset{r}{T}-\overset{r-1}{T}\sim \overset{r+1}{T}-\overset{r}{T}, \\ 
& r=1,\dots,k-1,\ T\to\infty. 
\end{split}
\ee 
\end{mydef81} 

\begin{mydef82}
Simultaneously with the Property 1, the sequence (\ref{2.8}) defines the partition of the integral 
\be \label{2.11} 
\int_T^{\overset{k}{T}}\left|\zf\right|^2{\rm d}t
\ee 
into the parts 
\be \label{2.12} 
\int_T^{\overset{k}{T}}\left|\zf\right|^2{\rm d}t=\sum_{r=1}^k\int_{\overset{r-1}{T}}^{\overset{r}{T}}\left|\zf\right|^2{\rm d}t, 
\ee 
that are asymptotically equal 
\be \label{2.13} 
\int_{\overset{r-1}{T}}^{\overset{r}{T}}\left|\zf\right|^2{\rm d}t\sim \int_{\overset{r}{T}}^{\overset{r+1}{T}}\left|\zf\right|^2{\rm d}t,\ T\to\infty. 
\ee 
\end{mydef82} 

It is clear, that (\ref{2.10}) follows from (\ref{2.3}) and (\ref{2.5}) since 
\be \label{2.14} 
\overset{r}{T}-\overset{r-1}{T}\sim (1-c)\frac{\overset{r}{T}}{\ln \overset{r}{T}}\sim (1-c)\frac{T}{\ln T},\ r=1,\dots,k, 
\ee  
while our eq. (\ref{2.13}) follows from (\ref{2.6}) and (\ref{2.5}).  

\section{Transformation of the basic formula (\ref{1.6})} 

\subsection{} 

Let us remind the assumption\footnote{See Remark 3.} 

\be \label{3.1} 
\{f_m(t)\}_{m=1}^\infty:\ f_m(t)\in C([0,2l]),\ m\in\mbb{N}. 
\ee 
Now, it follows from (\ref{1.6}) by mean value theorem that 
\be \label{3.2} 
\begin{split}
& \frac{1}{\prod_{r=0}^{k-1}\omega[\vp_1^r(\alpha)]}\int_{\overset{k}{T}}^{\overset{k}{\wideparen{T+2l}}}f_m^2(\vp_1^k(t)-T)\prod_{r=0}^{k-1}\left|\zeta\left(\frac 12+i\vp_1^r(t)\right)\right|^2{\rm d}t= \\ 
& A_m,\ m\in\mbb{N}, \\ 
& \omega[\vp_1^r(t)]=\left\{1+\mcal{O}\left(\frac{\ln\ln\vp_1^r(\alpha)}{\ln\vp_1^r\alpha}\right)\right\}\ln\vp_1^r(\alpha), 
\end{split}
\ee 
where 
\be \label{3.3} 
\alpha\in[\overset{k}{T},\overset{k}{\wideparen{T+2l}}],\ \alpha=\alpha(\overset{k}{T},\overset{k}{\wideparen{T+2l}},[\zeta],[\vp_1^0],\dots,[\vp_1^k]). 
\ee 

The following statement holds true: If 
\be \label{3.4} 
\alpha\in[\overset{k}{T},\overset{k}{\wideparen{T+2l}}], 
\ee  
then 
\be \label{3.5} 
\begin{split}
& \vp_1^0(\alpha)=\alpha\in [\overset{k}{T},\overset{k}{\wideparen{T+2l}}], \\ 
& \vp_1^1(\alpha)\in [\overset{k-1}{T},\overset{k-1}{\wideparen{T+2l}}], \\ 
& \vp_1^2(\alpha)\in [\overset{k-2}{T},\overset{k-2}{\wideparen{T+2l}}], \\ 
& \vdots \\ 
& \vp_1^{k-1}(\alpha)\in [\overset{1}{T},\overset{1}{\wideparen{T+2l}}], 
\end{split}
\ee  
for every fixed $l>0$, $l$ is independent on $T$, see Remark 2, i.e. 
\be \label{3.6} 
l=o\left(\frac{T}{\ln T}\right). 
\ee  

Next, let us remind that\footnote{See (\ref{2.4}), (\ref{2.14}).} 
\be \label{3.7} 
\overset{1}{T}<\overset{2}{T}<\dots<\overset{k}{T}, 
\ee  
and 
\be \label{3.8} 
\begin{split}
& \overset{r}{T}-\overset{r-1}{T}\sim (1-c)\frac{\overset{r}{T}}{\ln \overset{r}{T}}\sim(1-c)\frac{T}{\ln T}, \\ 
& r=1,\dots,k+1. 
\end{split}
\ee  
Then, by summation in (\ref{3.8}), we obtain 
\be \label{3.9} 
\overset{k+1}{T}-T\sim (k+1)(1-c)\frac{T}{\ln T}=\mcal{O}\left(k\frac{T}{\ln T}\right)\equiv \mcal{O}_k\left(\frac{T}{\ln T}\right), 
\ee  
i.e. 
\be \label{3.10} 
\overset{r}{T}-T\sim \mcal{O}_k\left(\frac{T}{\ln T}\right),\ r=1,\dots,k+1. 
\ee 
Now, it is true that\footnote{See (\ref{3.5}), (\ref{3.6}) and (\ref{3.10}).} 
\be \label{3.11} 
\begin{split}
& \vp_1^{k-1}(\alpha)\in[\overset{1}{T},\overset{2}{T}], \\ 
& \vp_1^{k-2}(\alpha)\in[\overset{2}{T},\overset{3}{T}], \\ 
& \vdots \\ 
& \vp_1^{0}(\alpha)=\alpha\in[\overset{k}{T},\overset{k+1}{T}], 
\end{split}
\ee 
and, of course, 
\be \label{3.12} 
\vp_1^r(\alpha)-T\leq \overset{k+1}{T}-T=\mcal{O}_k\left(\frac{T}{\ln T}\right),\ r=1,\dots,k+1. 
\ee 
Next 
\be \label{3.13} 
\vp_1^r(\alpha)=T+\{\vp_1^r(\alpha)-T\}=T\let\{1+\frac{\vp_1^r(\alpha)-T}{T}\}, 
\ee  
and\footnote{See (\ref{3.12}).} 
\be \label{3.14} 
\begin{split}
& \ln\vp_1^r(\alpha)=\ln T+\mcal{O}_k\left(\frac{1}{\ln T}\right)=\left\{1+\mcal{O}_k\left(\frac{1}{\ln^k T}\right)\right\}\ln T, \\ 
& \ln\ln\vp_1^r(\alpha)=\ln\ln T+\mcal{O}_k\left(\frac{1}{\ln^2 T}\right), 
\end{split} 
\ee 
i.e. 
\be \label{3.15} 
\frac{\ln\ln\vp_1^r(\alpha)}{\ln\vp_1^r(\alpha)}=\mcal{O}_k\left(\frac{\ln \ln T}{\ln T}\right),\ r=1,\dots,k+1. 
\ee 
Consequently\footnote{See (\ref{3.2}), (\ref{3.14}) and (\ref{3.15}).},
\be \label{3.16} 
\begin{split}
& \prod_{r=0}^{k-1}\omega[\vp_1^r(\alpha)]=\left\{\prod_{r=0}^{k-1}\left[1+\mcal{O}_k\left(\frac{\ln\ln T}{\ln T}\right)\right]\right\}\ln^kT= \\ 
& \left\{1+\mcal{O}_{k^2}\left(\frac{\ln\ln T}{\ln T}\right)\right\}\ln^kT. 
\end{split}
\ee 
That means, we have the following result. 

\begin{mydef51}
By assumption (\ref{3.1}) it is true that 
\be \label{3.17} 
\prod_{r=0}^{k-1}\omega[\vp_1^r(\alpha)]=\left\{1+\mcal{O}_{k^2}\left(\frac{\ln\ln T}{\ln T}\right)\right\}\ln^kT
\ee 
for every fixed $k\in\mbb{N}$. 
\end{mydef51}  
Next lemma follows immediately from (\ref{3.2}) by (\ref{3.17}). 
\begin{mydef52}
By assumption (\ref{3.1}) it is true that 
\be \label{3.18} 
\begin{split}
& \int_{\overset{k}{T}}^{\overset{k}{\wideparen{T+2l}}}f_m^2(\vp_1^k(t)-T)\prod_{r=0}^{k-1}\left|\zeta\left(\frac 12+i\vp_1^r(t)\right)\right|^2{\rm d}t= \\ 
& \left\{1+\mcal{O}_{k^2}\left(\frac{\ln\ln T}{\ln T}\right)\right\}A_m\ln^kT
\end{split}
\ee 
for every fixed $l>0$, $k,m\in\mbb{N}$. 
\end{mydef52}  

\section{New functionals and corresponding new $\zeta$-equivalents of the Fermat-Wiles theorem generated by the basic formula (\ref{1.6})}

\subsection{} 

We use the substitution 
\be \label{4.1} 
T=\exp\left\{\frac{x}{\sqrt{A_m}}\ln\tau\right\} 
\ee  
and the symbol 
\bdis 
W(x,\tau)=W(x,\tau;k,m)=\exp\left\{\frac{x}{\sqrt[k]{A_m}}\ln\tau\right\}
\edis 
in the eq. (\ref{3.18}). That gives the following functional. 

\begin{mydef11}
It is true by the assumption (\ref{3.1}) that 
\be \label{4.2} 
\begin{split}
& \lim_{\tau\to\infty}\frac{1}{\ln^k\tau}\int_{[W(x,\tau)]^k}^{[W(x,\tau)+2l]^k}f_m^2(\vp_1^k(t)-W(x,\tau))\times \\ 
& \prod_{r=0}^{k-1}\left|\zeta\left(\frac 12+i\vp_1^r(t)\right)\right|^2{\rm d}t=x^k, 
\end{split}
\ee  
for every fixed 
\bdis 
x>0,\ l>0,\ k,m\in\mbb{N}. 
\edis 
\end{mydef11}

Next, the substitution 
\bdis 
x\to\FR
\edis 
into (\ref{4.2}) gives the following result. 

\begin{mydef12}
The $\zeta$-condition 
\be \label{4.3} 
\begin{split}
& \lim_{\tau\to\infty}\frac{1}{\ln^k\tau}\int_{[W(\FR,\tau)]^k}^{[W(\FR,\tau)+2l]^k}f_m^2(\vp_1^k(t)-W(\FR,\tau))\times \\ 
& \prod_{r=0}^{k-1}\left|\zeta\left(\frac 12+i\vp_1^r(t)\right)\right|^2{\rm d}t\not=1, 
\end{split}
\ee 
on the set of all Fermat's rationals expresses the next $\zeta$-equivalent of the Fermat-Wiles theorem for every fixed $l>0$ and $k,m\in\mbb{N}$. 
\end{mydef12} 

\section{Fuctionals and $\zeta$-equivalents of the Fermat-Wiles theorem that follow from the elementary Fourier orthogonal system} 

\subsection{} 

Lemma 1 is applicable also in the cases (\ref{1.13}) -- (\ref{1.15}) and this leads to the following results. 

\begin{itemize}
	\item[(A)] The formula (\ref{1.13}) generates the functional: 
	\begin{mydef13}
	\be \label{5.1} 
	\lim_{T\to\infty}\frac{1}{\ln^kT}\int_{\overset{k}{T}}^{\overset{k}{\wideparen{T+2l}}}\prod_{r=0}^{k-1}\left|\zeta\left(\frac 12+i\vp_1^r(t)\right)\right|^2{\rm d}t=2l, 
	\ee  
	for every fixed $l>0$ and $k\in\mbb{N}$. 
	\end{mydef13} 
	\begin{remark}
	The functional $\mcal{F}_1$ is defined by the formula (\ref{5.1}) as follows: 
	\be \label{5.2} 
	T+2l\xrightarrow{\mcal{F}_1}2l. 
	\ee 
	\end{remark}
	Next, in the case\footnote{See (\ref{1.8}).}
	\bdis 
	l\to\FR 
	\edis  
	we have the following result. 
	\begin{mydef14}
	The $\zeta$-condition 
	\be \label{5.3} 
	\lim_{T\to\infty}\frac{1}{\ln^kT}\int_{\overset{k}{T}}^{[T+2\FR]^k}\prod_{r=0}^{k-1}\left|\zeta\left(\frac 12+i\vp_1^r(t)\right)\right|^2{\rm d}t\not=2
	\ee 
	on the set of all Fermat's rationals expresses the next $\zeta$-equivalent of the Fermat-Wiles theorem for every fixed $k\in\mbb{N}$. 
	\end{mydef14} 
	\item[(B)] The formula (\ref{1.14}) generates the following functional. 
	\begin{mydef15} 
	\be \label{5.4} 
	\begin{split}
	& \lim_{T\to\infty}\frac{1}{\ln^kT}\int_{\overset{k}{T}}^{\overset{k}{\wideparen{T+2l}}}\prod_{r=0}^{k-1}\left|\zeta\left(\frac 12+i\vp_1^r(t)\right)\right|^2\times \\ 
	& \cos^2(\frac{\pi}{l}m(\vp_1^k(t)-T)){\rm d}t=l 
	\end{split}
	\ee 
	for every fixed $l>0$ and $k,m\in\mbb{N}$. 
	\end{mydef15} 
	\begin{remark}
	The functional $\mcal{F}_2$ is defined by the formula (\ref{5.4}) as follows: 
	\be \label{5.5} 
	T+2l\xrightarrow{\mcal{F}_2}l. 
	\ee
	\end{remark} 
	Next, in the case of Fermat's rationals, we have the following theorem.
	\begin{mydef16}
	The $\zeta$-condition 
	\be \label{5.6} 
	\begin{split}
		& \lim_{T\to\infty}\frac{1}{\ln^kT}\int_{\overset{k}{T}}^{[T+2\FR]^k}\prod_{r=0}^{k-1}\left|\zeta\left(\frac 12+i\vp_1^r(t)\right)\right|^2\times \\ 
		& \cos^2(\pi m\FRI(\vp_1^k(t)-T)){\rm d}t\not=1 
	\end{split}
	\ee 
	on the set of all Fermat's rationals expresses the next $\zeta$-equivalent of the Fermat-Wiles theorem for every fixed $k,m\in\mbb{N}$. 
	\end{mydef16}
\end{itemize} 

\subsection{} 

We give here some subsidiary remarks. 

\begin{remark}
The differents of the two limits (\ref{5.1}) and (\ref{5.4}) gives, of course, the functional 
\be \label{5.7} 
	\begin{split}
	& \lim_{T\to\infty}\frac{1}{\ln^kT}\int_{\overset{k}{T}}^{\overset{k}{\wideparen{T+2l}}}\prod_{r=0}^{k-1}\left|\zeta\left(\frac 12+i\vp_1^r(t)\right)\right|^2\times \\ 
	& \sin^2(\frac{\pi}{l}m(\vp_1^k(t)-T)){\rm d}t=l. 
\end{split}
\ee 
\end{remark}

\begin{remark}
The difference of the two limits (\ref{5.4}) and (\ref{5.7}) gives us the following 
\be \label{5.8} 
\begin{split}
	& \lim_{T\to\infty}\frac{1}{\ln^kT}\int_{\overset{k}{T}}^{\overset{k}{\wideparen{T+2l}}}\prod_{r=0}^{k-1}\left|\zeta\left(\frac 12+i\vp_1^r(t)\right)\right|^2\times \\ 
	& \cos(\frac{\pi}{l}2m(\vp_1^k(t)-T)){\rm d}t=0,  
\end{split}
\ee 
since, of course, the integral can be reduced to the condition 
\be \label{5.9} 
1\perp\cos(\frac{\pi}{l}2mt) 
\ee 
on the segment $[0,2l]$, comp. (\ref{1.3}) and (\ref{1.11}). 
\end{remark} 

\section{New expression for $\ln^kT$ and some Remarks} 

\subsection{} 

Let us remind that the eq. (\ref{3.18}) gives, in the case 
\bdis 
f_m(t)=1,\ t\in [0,2l],\ l=\frac 12, 
\edis  
the following formula. 

\begin{mydef53}
\be \label{6.1} 
\int_{\overset{k}{T}}^{\overset{k}{\wideparen{T+2l}}}\prod_{r=0}^{k-1}\left|\zeta\left(\frac 12+i\vp_1^r(t)\right)\right|^2{\rm d}t=
\left\{1+\mcal{O}_{k^2}\left(\frac{\ln\ln T}{\ln T}\right)\right\}\ln^kT, 
\ee  
for every fixed $k\in\mbb{N}$. 
\end{mydef53} 

\begin{remark}
\be \label{6.2} 
\begin{split}
& \left\{\int_{\overset{k}{T}}^{\overset{k}{\wideparen{T+2l}}}\prod_{r=0}^{k-1}\left|\zf\right|^2{\rm d}\right\}^{1/k}=
\left\{1+\mcal{O}_{k}\left(\frac{\ln\ln T}{\ln T}\right)\right\}\ln T, 
\end{split}
\ee 
where\footnote{See (\ref{3.16}).} 
\bdis 
\mcal{O}_{k^2}\left(\frac{\ln\ln T}{\ln T}\right)=\mcal{O}\left(k^2\frac{\ln\ln T}{\ln T}\right). 
\edis 
\end{remark} 

\subsection{} 

The formula (\ref{6.2}) is a supplement to the set of our quotient formulae: 
\be \label{6.3} 
\frac{\int_T^{\overset{1}{T}}|\zf|^2{\rm d}t}{\int_T^{\overset{1}{T}}|\zeta(\sigma+it)|^2{\rm d}t}=\frac{1}{\zeta(2\sigma)}\ln T+\mcal{O}(1),\ \sigma\geq\frac 12+\epsilon, 
\ee 
and 
\be \label{6.4} 
\frac{\int_T^{\overset{1}{T}}|\zf|^2{\rm d}t}{\int_T^{\overset{1}{T}}|S_1(t)|^{2p}{\rm d}t}=\frac{1}{\bar{c}(p)}\ln T+\mcal{O}(1), 
\ee 
where $\epsilon>0$ is fixed small number, and\footnote{See \cite{14}, (3.11), \cite{17}, (4.5).} 
\be \label{6.5} 
S_1(t)=\frac{1}{\pi}\int_0^t\arg\zfu{\rm d}u. 
\ee 

\subsection{} 

What is our reason for the study of infinite processes in mathematics? An adequate answer is given by H. Weyl in his book \cite{18}. 
\begin{remark}
Mathematics is the science of infinite, its goal the symbolic comprehension of the infinite with human, that is finite, means. (\cite{18}, p. 7.) 
\end{remark} 

\begin{remark}
Indeed, God as completed infinite cannot and will not be comprehended by it; neither can God penetrate into man by revelation, nor man penetrate to him by mystical perception. The completed infinite we can only represent in symbols. From this relationship every creative act of man receives its deep consecration and dignity. But only in mathematics and physics, as far as I can see, has symbolical-theoretical construction acquired sufficient solidity to be convincing for everyone whose mind is open to these sciences. (\cite{18}, p. 84.)
\end{remark}

I would like to thank Michal Demetrian for his moral support of my study of Jacob's ladders.

\end{document}